\documentclass[pamm,a4paper,fleqn]{w-art}
\usepackage{times,cite,w-thm}
\usepackage[T1]{fontenc}
\usepackage[utf8]{inputenc}
\usepackage{comment}
\theoremstyle{plain}

\theoremstyle{definition}

\usepackage{graphicx, subcaption}
\usepackage{subfig}
\usepackage{amssymb}
\usepackage{amsfonts}
\usepackage{amsmath}
\usepackage{algorithm}
\usepackage{algpseudocode}
\begin{document}

\TitleLanguage[EN]
\title[The short title]{Optimization of Approximate Maps for Linear Systems Arising in Discretized PDEs}


\author{\firstname{Rishad}  \lastname{Islam}\inst{1,}%
  \footnote{Corresponding author: mds222@lehigh.edu}}

\address[\inst{1}]{\CountryCode[DE]Lehigh University, Bethlehem, Pennsylvania, USA}
\author{\firstname{Arielle} \lastname{Carr}\inst{1,}%
  \footnote{arg318@lehigh.edu}}

\author{\firstname{Colin} \lastname{Jacobs}\inst{1,}%
  \footnote{cbj225@lehigh.edu}}
\AbstractLanguage[EN]
\begin{abstract}
Generally, discretization of partial differential equations (PDEs) creates a sequence
of linear systems $A_k x_k = b_k, k = 0, 1, 2, ..., N$ with well-known and structured sparsity
patterns. Preconditioners are often necessary to achieve fast convergence When solving these linear 
systems using iterative solvers. We can use preconditioner updates for closely related systems instead
of computing a preconditioner for each system from scratch. One such preconditioner update is the
sparse approximate map (SAM), which is based on the sparse approximate inverse
preconditioner using a least squares approximation. A SAM then acts as a map from one matrix
in the sequence to another nearby one for which we have an effective preconditioner. To efficiently
compute an effective SAM update (i.e., one that facilitates fast convergence of the iterative solver),
we seek to compute an optimal sparsity pattern. In this paper, we examine several sparsity patterns
for computing the SAM update to characterize optimal or near-optimal sparsity patterns
for linear systems arising from discretized PDEs. 
\end{abstract}
\maketitle                   

\section{Introduction}
Methods such as finite elements or finite differences are used to solve partial differential equations (PDE) numerically. These discretization processes often generate a sequence of linear systems, $A_kx_k=b_k$, $k = 0,1,2,\dots N$ that require solving to obtain the solution to the underlying PDE. Iterative methods, such as Krylov subspace methods\cite{Vand03, Saad03} are effective techniques for solving these linear systems. However, the systems generally must be well-conditioned to achieve fast convergence of the Krylov methods, so we often employ mitigation strategies such as preconditioning to reduce the condition number of the coefficient matrix (see, e.g., \cite{benzi2002preconditioning} for an overview of preconditioning techniques). Unfortunately, computing and applying an effective preconditioner can be very expensive, and in the case of a long sequence of ill-conditioned linear systems, doing so for each system quickly becomes computationally infeasible.

Sparse Approximate Maps (SAMs) \cite{carr2021preconditioning} allow us to circumvent some of this cost by acting as a preconditioner update and can be used irrespective of the initial preconditioner type. They are shown to be particularly effective for expensive preconditioners like the algebraic multigrid preconditioner. The computation of the SAM update is based on the Sparse Approximate Inverse (SAI) \cite{benzi1999comparative} and acts as an approximate map between two matrices in a sequence: from one later in the sequence to another for which we have previously computed an effective preconditioner. The efficiency of a SAM depends not only on its efficacy as a preconditioner but also on the apriori choice of the sparsity patterns of the maps\cite{carr2021preconditioning}. However, while several different sparsity patterns are empirically evaluated in \cite{carr2021preconditioning}, only a fixed sparsity pattern is ultimately selected as the best choice for the maps throughout the entire sequence. For changing matrices in the sequence, this might not adequately represent the evolving behavior of the systems. In particular, the sparsity pattern should account for the changing matrices in the sequence, both in terms of the number and location of nonzero elements as well as the magnitude of the matrix entries. 

In this paper, we analyze different choices of sparsity patterns for SAMs corresponding to each matrix in the sequence of linear systems to more carefully characterize optimal choices in sparsity patterns for approximate maps between linear systems arising in discretized PDEs. We consider the choice of sparsity patterns for the SAM updates and how those choices affect the residual of the approximate map. We examine the graphical representation of the maps to understand their behavior and select our choices for sparsity patterns based on these observations. We also study the corresponding sparsity patterns of the maps for the changing matrices in the sequence and compare them with the fixed sparsity pattern used in the previous version of SAM in \cite{carr2021preconditioning}. Extending the analysis in \cite{chow2000priori}, we show that the sparsity pattern of the exact map is a subset of the sparsity pattern of the transitive closure of a graph representation of $A_k$, denoted $G(A_k)$. As a result, we focus our attention on level 1 and level 2 neighbors of $G(A_k)$, as well as several sparsification approaches based directly on the sparsity pattern of $A_k$. We modify the current implementation of the SAM update to account for the changing sparsity patterns and matrix entries of $A_k$ in the sequence and show that if our sequence of linear systems satisfies simple subset assumptions, we can characterize properties of the residual of the SAM.

\section{Sparse Approximate Maps}
The main idea behind the SAM update is to approximately map one coefficient matrix, $A_\ell$, in a sequence of linear systems to another matrix, $A_j$ for $j < \ell$, for which we already have an effective preconditioner. This allows us to reuse the preconditioner computed for $A_j$ at a dramatically reduced cost (see \cite{carr2021preconditioning} for an extensive experimental analysis using several PDE-based linear systems).  In this paper, we consider the specific case of mapping each $k^{th}$ matrix in the sequence back to the initial matrix $A_0$. In our subsequent analysis of the quality of the map, we define $A_k$ as the source matrix and $A_0$ as the target matrix. We also assume $A_k \in \mathbb{C}^{n \times n}$ and $x_k, b_k \in \mathbb{C}^n$.

First, we define the exact map between matrices as provided in \cite{carr2021preconditioning}: given a sequence of matrices $A_0, A_1, A_2, \dots,N$, the exact map $\widehat{N_k}$ for each system is
\begin{equation}
  \label{exact_map}
  A_k\widehat{N_k}=A_0.
\end{equation}
Computing (\ref{exact_map}) generally results in a dense matrix that is far too computationally expensive to use in practice. Therefore, we instead compute an approximate map $N_k$, or SAM update, to reduce the computational costs associated with determining $\widehat{N_k}$ while maintaining optimal properties when mapping $A_k$ back to $A_0$. Formally, the SAM update, $N_k$, for the $k^{th}$ linear system, $A_k x_k = b_k$, is defined as
\begin{equation}
  \label{approx_map}
  N_k = \arg \min_N ||A_kN - A_0||_F
\end{equation}
where $N_k$ satisfies an imposed sparsity pattern $S$ for $\mathbb{C}^{n\times n}$. $S$ can be defined as any subset of $\{1,2,\dots,n\}\times \{1,2,\dots,n\}$. Throughout this paper, we use the notation $S(B)$, to denote ``the sparsity pattern of the matrix $B$." From (\ref{approx_map}), we also define the residual of the SAM as 
\begin{equation}
  \label{res}
  R_k = A_kN_k - A_0.
\end{equation}

In fact, (\ref{approx_map}) results in $n$ small least squares problems, each of which can be solved by choosing a smaller submatrix from $A_k$ based on the sparsity pattern $S$ for each column of $N_k$ and solving the smaller linear system directly, where the RHS is the corresponding column of $A_0$. In other words, we solve $n_j = \min_m \|\mathcal{A}_k m - (A_0)_j\|_F$ for $j = 1, 2, \dots, n$, where $\mathcal{A}$ represents the submatrix of $A_k$ subject to the chosen sparsity pattern, $(A_0)_j$ is the $j^{th}$ column of $A_0$, and $n_j$ is the resulting $j^{th}$ column of the SAM, or approximate map. So, the effectiveness of this approximate map depends on determining optimal sparsity patterns\cite{carr2021preconditioning}, whereby we must optimize the number and location of the nonzeros to optimize the accuracy of the map without negatively affecting the cost to compute it. The sparsity pattern should be sparse enough to enable fast computation yet contain the appropriate amount and location of nonzeros to provide an accurate approximation to the exact map $\widehat{N_k}$. In general, adaptive sparsity patterns are more accurate but tend to be more computationally expensive\cite{benzi2002preconditioning}. Therefore, as in \cite{carr2021preconditioning, chow2000priori}, we focus our attention on a priori sparsity patterns in the current paper. However, future work will seek to extend investigations in, e.g., \cite{ChowSaad98, GrotHuck97} to improve the efficiency of adaptive strategies, particularly in high-performance computing environments by taking advantage of the inherently parallel nature of the SAM computation. In Section \ref{sec:optimalSp}, we examine several sparsification strategies based on the sparsity pattern of $A_k$ but first motivate our choice in strategies based on a graph interpretation of the exact map $\widehat{N_k}$.

\section{Graphical Representation of Exact Maps}\label{sec:graph}
In this section, we will observe the properties of the exact map through the adjacency graph representation of the matrices, specifically that of the source matrix, $A_k$. In the current implementation of the algorithm \cite{carr2021preconditioning}, the sparsity patterns used are based on that of the target matrix, $A_0$. We theorize that if the matrices in the sequence change according to a specific property (e.g., in response to a physical change that the underlying PDE is describing), it may be more judicious to use a sparsity pattern based on that of the source matrix. 

We motivate our analysis using a graph interpretation of the exact map based on the analysis of sparsity patterns for SAIs in \cite{chow2000priori} and adopt the notation used therein. The adjacency graph representation of a matrix $A$ of order $n$ is defined as the directed graph $G(A)$ with vertices $1, 2, \dots, n$ and whose edges $(i \rightarrow j)$ correspond to non-zero off-diagonal entries in $A$. A subset of $G(A)$ contains the same vertices as $G(A)$ and a subset of edges from $G(A)$. Thus, the graph representation of the nonzeros in column $j$ of matrix $A$ is the set of vertices in $G(A)$ with edges pointing to vertex $j$, including the vertex $j$ itself. This set of vertices can be considered the level 1 neighbor of column $j$. The level 2 neighbors of column $j$ will be the set of vertices in $G(A)$ with a path of length at most 2 pointing to vertex $j$. As discussed in cite{chow2000priori}, the graph representation of the inverse of a matrix $A^{-1}$ contains all the possible edges between all the vertices in $G(A)$, known as the transitive closure of the matrix $A$ and denoted as $G^*(A)$. This can also be represented by $G^(A^n)$. It is important to highlight that this does not imply $A^{-1} = A^n$, but rather that the {\it sparsity pattern} of $A^{-1}$ is the same as that of $A^n$, which can then be distilled from $G^*(A)$. In \cite{chow2000priori}, this conclusion is used to justify using sparsity patterns when computing the SAI that are based on subsets of the sparsity pattern of $A$. This becomes immediately relevant to our analysis of sparsity patterns for SAMs. 

First, we assume that the sparsity pattern of the sequence of matrices satisfies
\begin{equation}
    \label{property}
    S(A_0) \subseteq S(A_1) \subseteq \dots \subseteq S(A_N),
\end{equation}
which is often the case for many applications involving discretized PDEs. From (\ref{exact_map}), we define the formula for the exact map as $\widehat{N_k} = A_k^{-1}A_0$. Note that the multiplication of two matrices $A$ and $B$ with the same set of vertices $1, 2, \dots, n$ can be represented by the graphical representation $G(AB)$. It contains all the paths between vertex $i$ and $j$ of $G(AB)$ using edges from $G(A)$ and $G(B)$ and can be interpreted as including any extra edges between the vertices in $G(A)$ according to the existing edges between the same vertices in $G(B)$. In the case where $S(A) \subseteq S(B)$, no new edges will be added to $G(AB)$. Thus, when considering the exact map $\widehat{N_k}$, with (\ref{property}) satisfied, it must follow that the multiplication of $A_k^{-1}$ by $A_0$ contributes no new directed edges to $G(\widehat{N_k})$.  In other words, for sequences of matrices that satisfy (\ref{property}), $G(\widehat{N_k}) = G^*(A_k)$, or $S(\widehat{N_k}) = S(A_k^n)$, and we can reasonably apply to SAMs the sparsification choices made in \cite{chow2000priori} and \cite{chow2001parallel} for SAIs.

\begin{figure}[h]
\begin{minipage}{80mm}
\captionsetup{justification=centering}
\includegraphics[width=\linewidth,height=50mm]{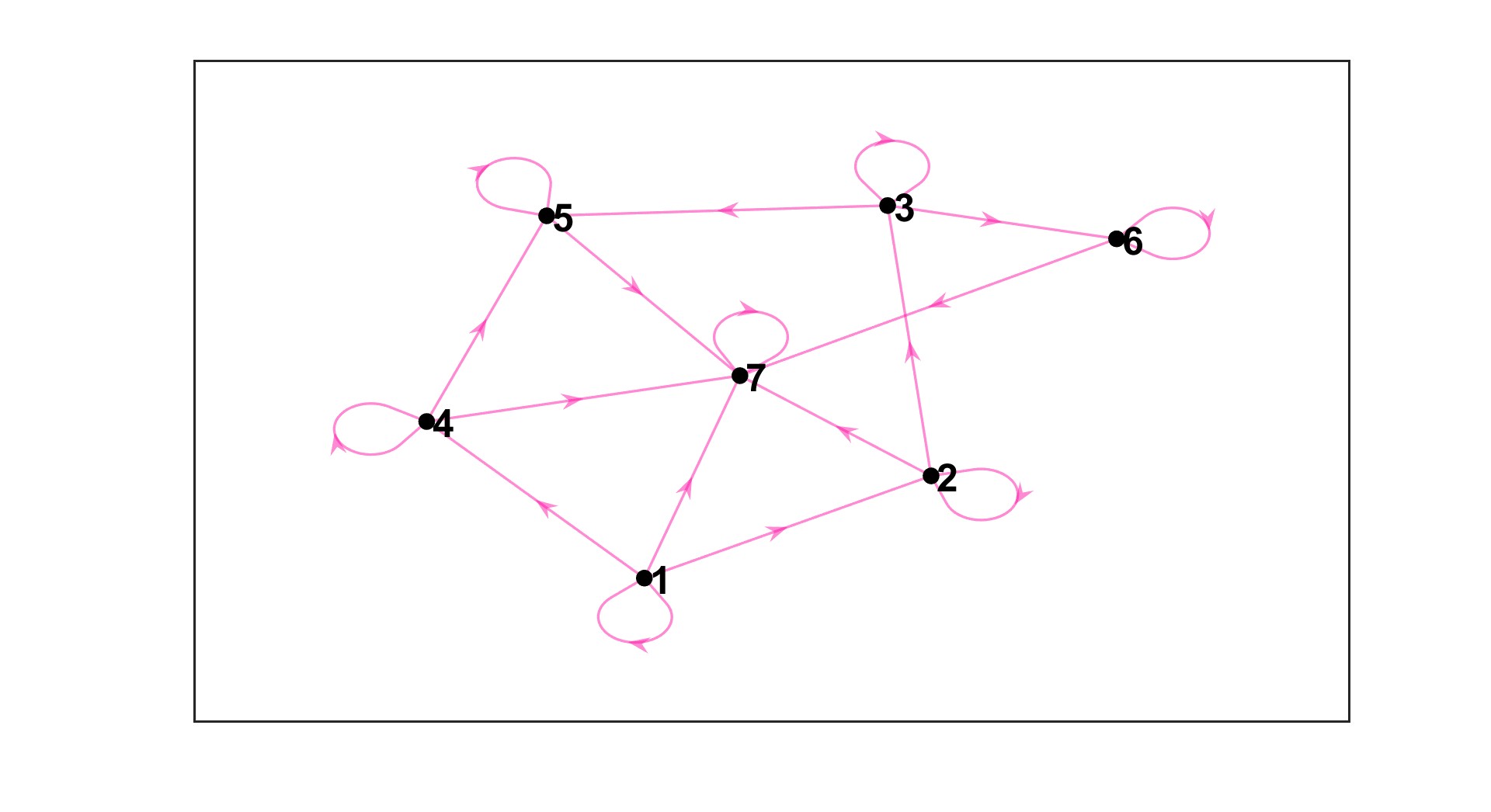}
\caption{Adjacency graph of target matrix, $G(A_0)$}
\label{G(A_0)}
\end{minipage}
\hfil
\begin{minipage}{80mm}
\captionsetup{justification=centering}
\includegraphics[width=\linewidth,height=50mm]{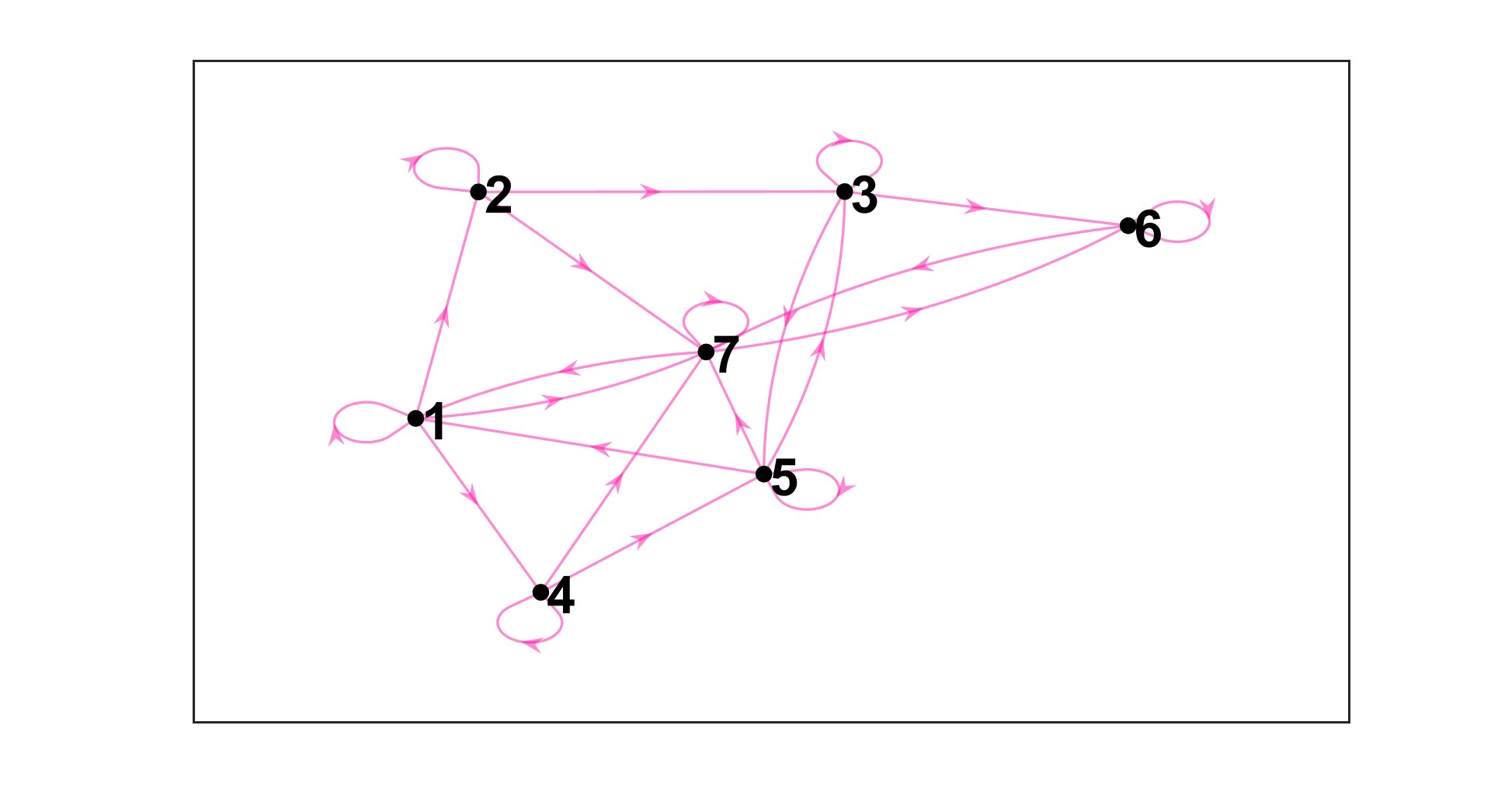}
\caption{Adjacency graph of source matrix, $G(A_1)$}
\label{G(A_k)}
\end{minipage}

\begin{minipage}{80mm}
\captionsetup{justification=centering}
\includegraphics[width=\linewidth,height=50mm]{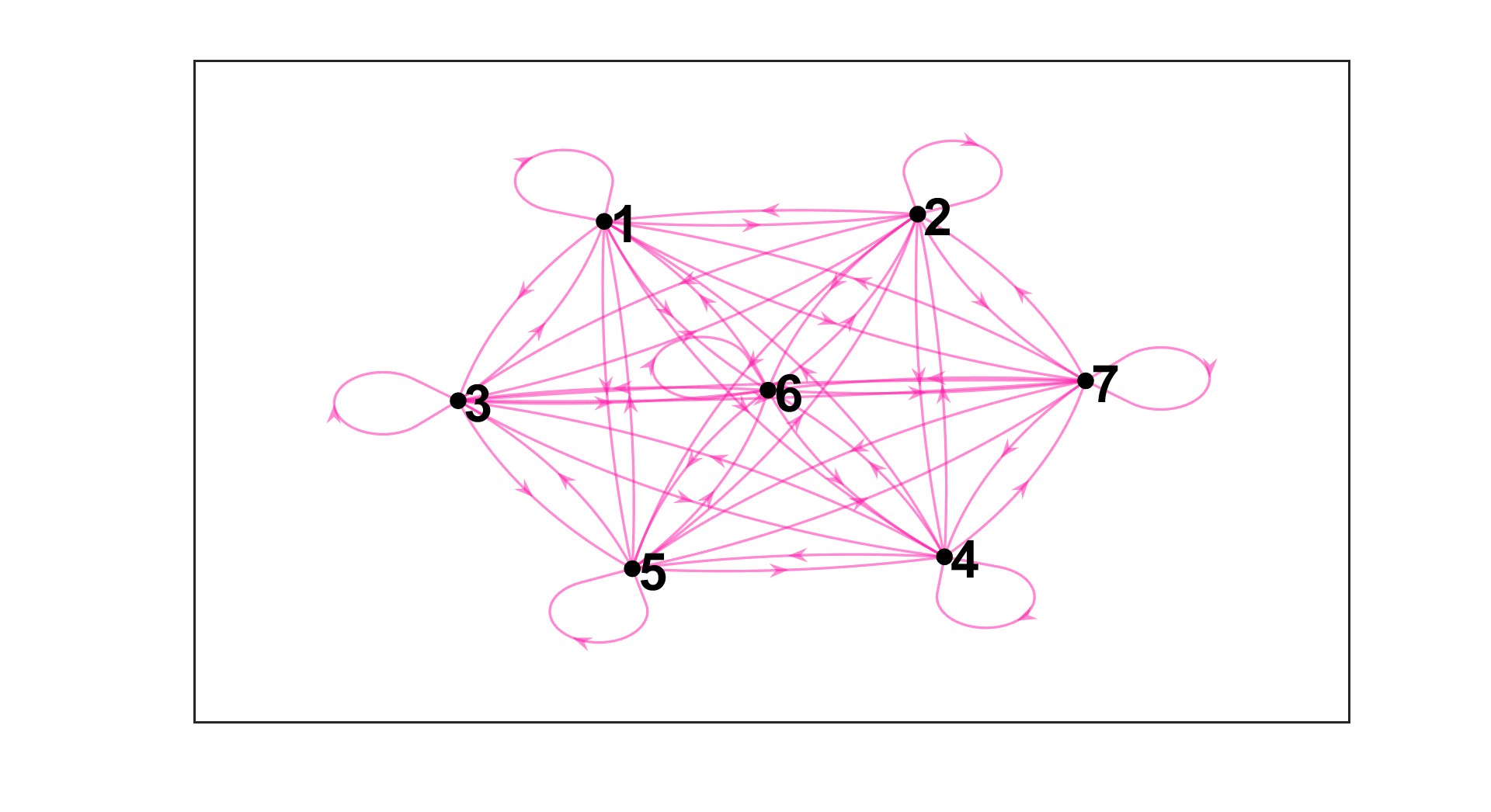}
\caption{Transitive closure of source matrix, $G^*(A_1)$}
\label{G(A_k_inv)}
\end{minipage}
\hfil
\begin{minipage}{80mm}
\captionsetup{justification=centering}
\includegraphics[width=\linewidth,height=50mm]{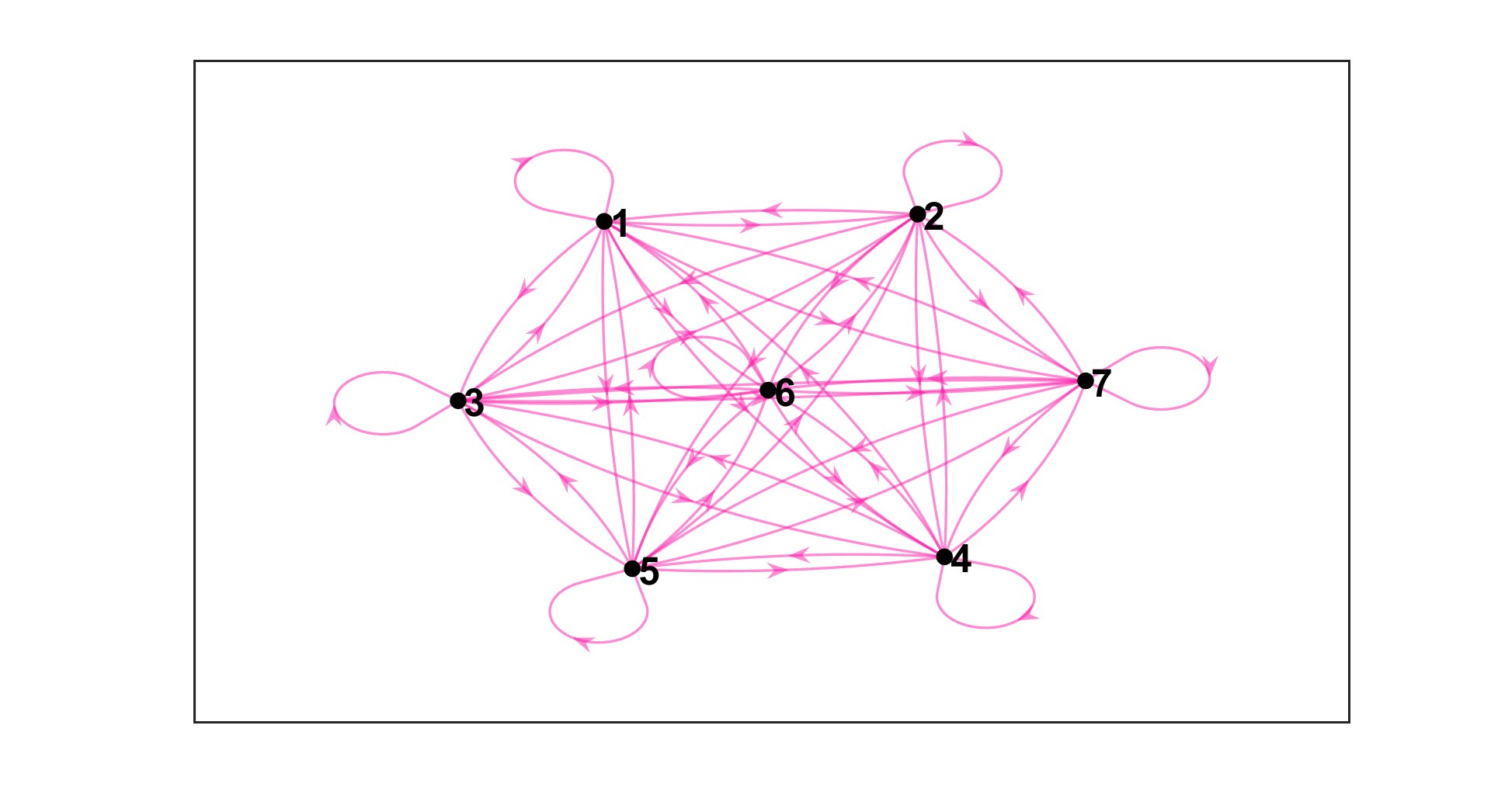}
\caption{Adjacency graph of exact map, $G(\widehat{N_1})$.}
\label{G(N_k)}
\end{minipage}
\end{figure}

We verify this conclusion using two example sparse matrices $A_0$ and $A_1$ with $n = 7$ and random non-zero entries such that $S(A_0) \subseteq S(A_1)$. Here, our analysis is not based on magnitude, so we display the boolean matrices in (\ref{mats}). 
\begin{equation}\label{mats}
A_0 = \begin{pmatrix}
1 & 1 & 0 & 1 & 0 & 0 & 1 \\
0 & 1 & 1 & 0 & 0 & 0 & 1 \\
0 & 0 & 1 & 0 & 1 & 1 & 0 \\
0 & 0 & 0 & 1 & 1 & 0 & 1 \\
0 & 0 & 0 & 0 & 1 & 0 & 1 \\
0 & 0 & 0 & 0 & 0 & 1 & 1 \\
0 & 0 & 0 & 0 & 0 & 0 & 1 
\end{pmatrix} \quad \quad \quad \quad
A_1 = \begin{pmatrix}
1 & 1 & 0 & 1 & 0 & 0 & 1 \\
0 & 1 & 1 & 0 & 0 & 0 & 1 \\
0 & 0 & 1 & 0 & 1 & 1 & 0 \\
0 & 0 & 0 & 1 & 1 & 0 & 1 \\
1 & 0 & 1 & 0 & 1 & 0 & 1 \\
0 & 0 & 0 & 0 & 0 & 1 & 1 \\
1 & 0 & 0 & 0 & 0 & 1 & 1 
\end{pmatrix}
\end{equation}
Figures \ref{G(A_0)} and \ref{G(A_k)} provide the graph representations of these two matrices. We can easily verify that the edges between the vertices in $G(A_0)$ are also present in $G(A_1)$. We explicitly compute the inverse of the source matrix $A_1$ and represent it using the transitive closure, $G^*(A_1)$ in figure \ref{G(A_k_inv)}. Similarly, we explicitly compute the exact map $\widehat{N_1} = A_1^{-1}A_0$ and represent it in figure \ref{G(N_k)}. Then, we compare the adjacency graph between the two figures \ref{G(A_k_inv)} and \ref{G(N_k)} and observe that they are identical. No new edges 
are added in $G(\widehat{N_1})$ after multiplying $A_1^{-1}$ with $A_0$. 

\section{Optimal Sparsity Patterns for Approximate Maps}\label{sec:optimalSp}
\begin{figure}[hh]
\begin{center}
\captionsetup{justification=centering}
\begin{subfigure}{.34\textwidth}
        		\includegraphics[width=\linewidth]{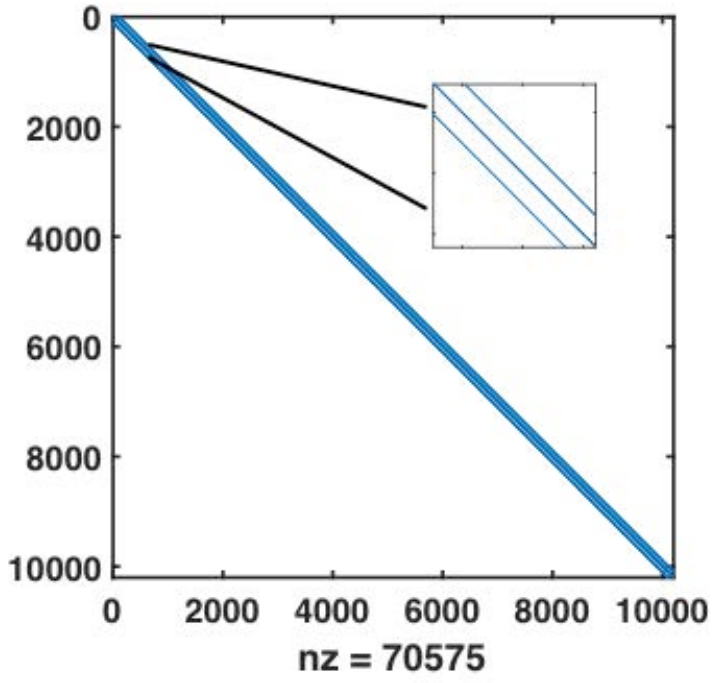}
		\subcaption{$S(A_k), k = 1, 2, \dots, 20$.}\label{fig:THTsparsity}
        \end{subfigure}
        
        \begin{subfigure}{.3\textwidth}
        		\includegraphics[width=\linewidth]{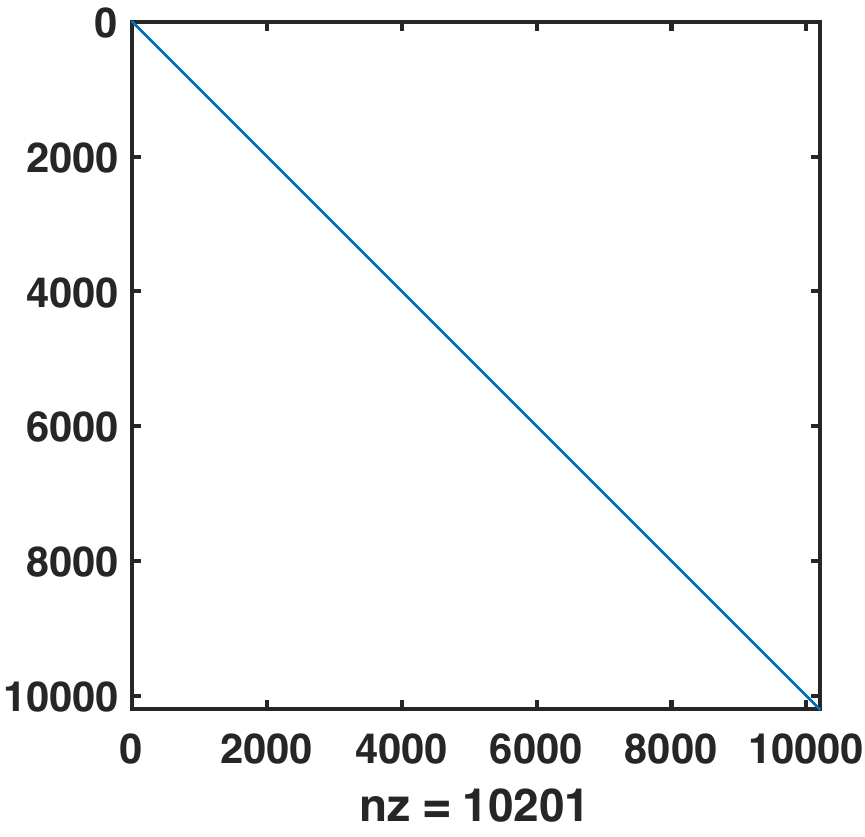}
		\subcaption{$\widehat{N_2} > 1.e-2$}
        \end{subfigure}
        \begin{subfigure}{.3\textwidth}
        		\includegraphics[width=\linewidth]{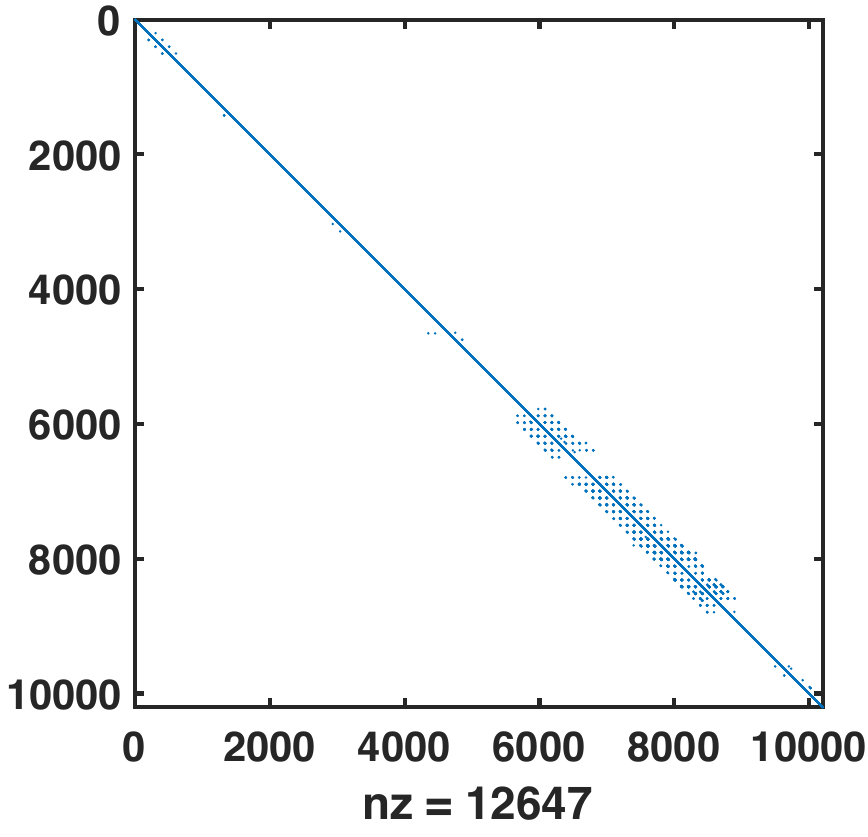}
		\subcaption{$\widehat{N_{10}} > 1.e-2$}
        \end{subfigure}
                \begin{subfigure}{.3\textwidth}
        		\includegraphics[width=\linewidth]{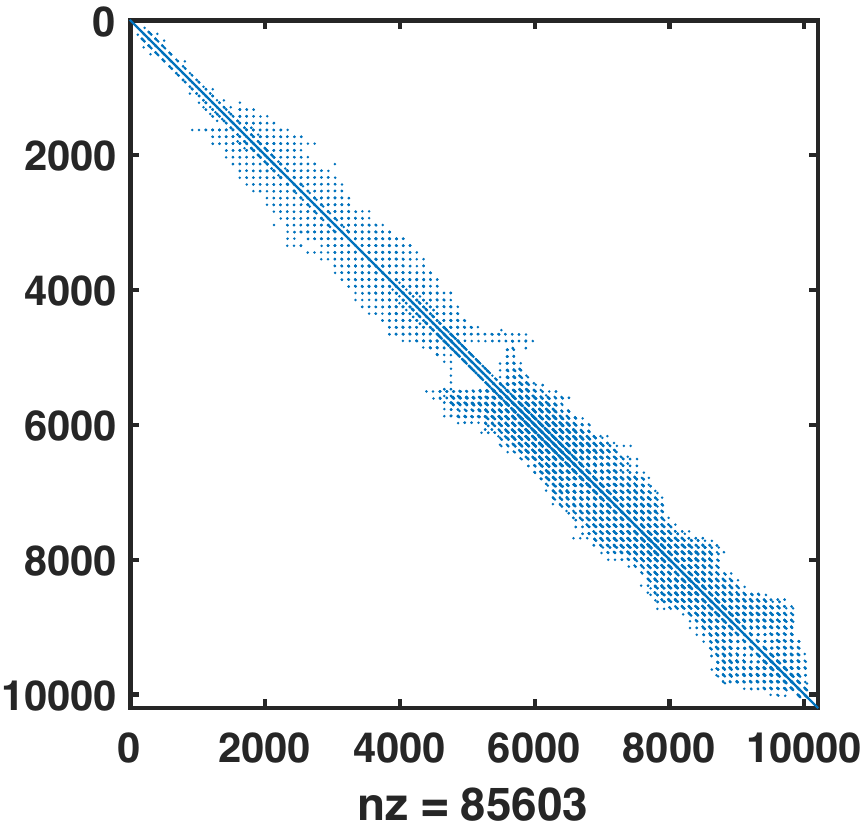}
		\subcaption{$\widehat{N_{20}} > 1.e-2$}
        \end{subfigure}

        \begin{subfigure}{.3\textwidth}
        		\includegraphics[width=\linewidth]{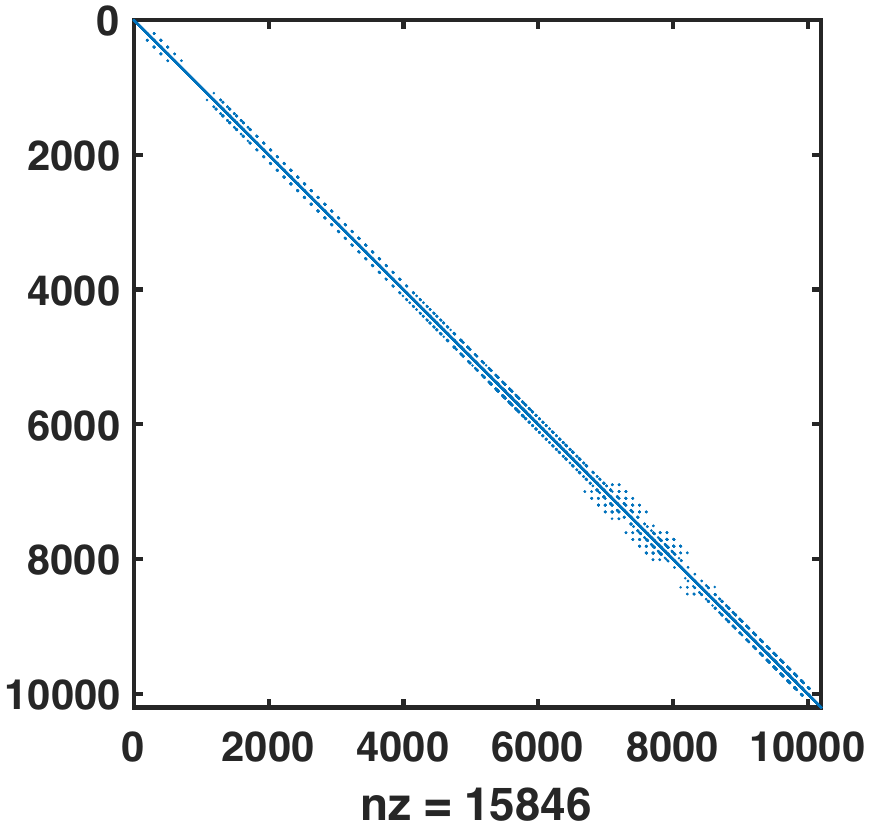}
		\subcaption{$\widehat{N_2} > 1.e-4$}
        \end{subfigure}
        \begin{subfigure}{.3\textwidth}
        		\includegraphics[width=\linewidth]{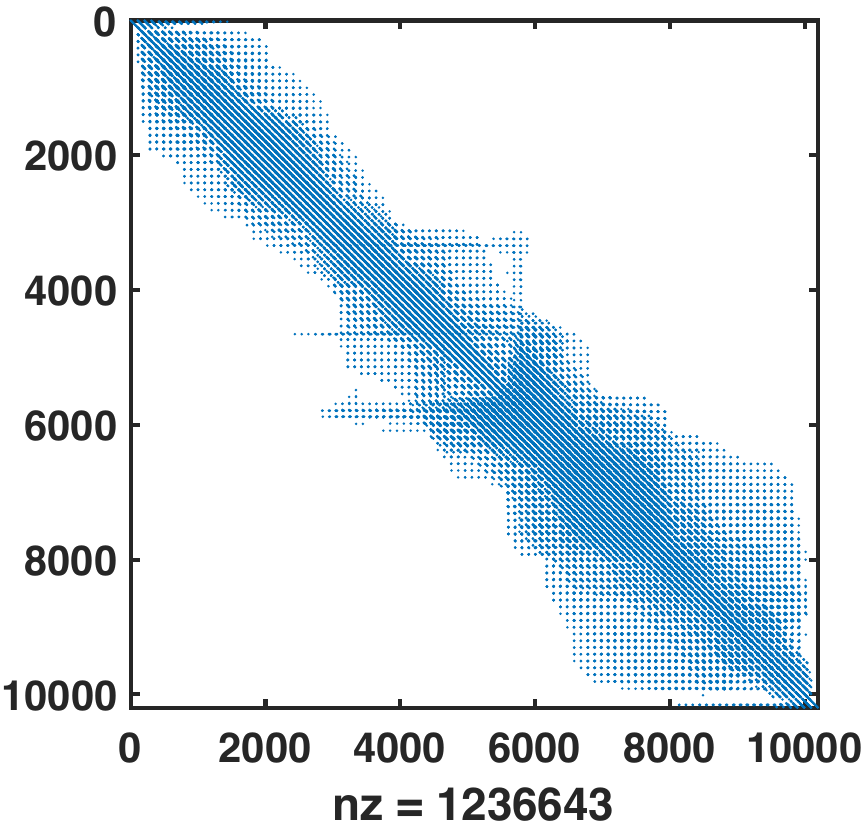}
		\subcaption{$\widehat{N_{10}} > 1.e-4$}
        \end{subfigure}
                \begin{subfigure}{.3\textwidth}
        		\includegraphics[width=\linewidth]{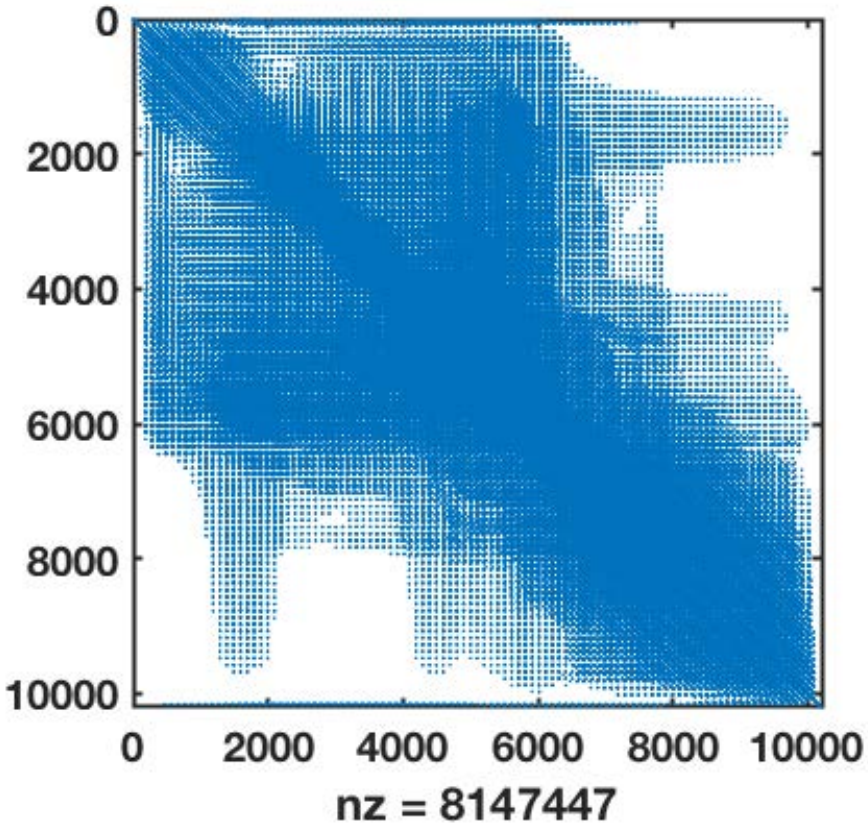}
		\subcaption{$\widehat{N_{20}} > 1.e-4$}
        \end{subfigure}
	  \caption{Sparsity patterns for THT matrices and global sparsification applied to the exact map,$\widehat{N_k}$, between selected matrices in the sequence. $\widehat{N_j} > 1.e-p$ denotes the dropping of elements (in absolute value) from the exact map that is less than $1.e-p$, for $p=2,4$.}
	  \label{THT}
\end{center}
\end{figure}

In addition to the graphical motivations discussed in Section \ref{sec}, we provide additional intuition when considering sparsification strategies applied directly to $S(\widehat{N_k})$. We use a sequence of matrices arising in transient hydraulic tomography (THT) applications \cite{cardiff20113}, a method for imaging the Earth's subsurface described by underlying PDEs that model groundwater flow. These PDEs are discretized via finite elements and result in a sequence of shifted systems that are solved at each time step. We show the resulting block tridiagonal sparsity pattern for these shifted systems in figure \ref{fig:THTsparsity}. For brevity, we refer the reader to \cite{Bakhos2015940, GriDesGug2020} for further details on this problem and a characterization of the specific PDEs defining this physical process. We consider three matrices, $A_0$, $A_{9}$, and $A_{19}$, extracted from a sequence of twenty THT matrices with $n = 10201$ that satisfy (\ref{property}) and compute the exact map for all the matrices in the sequence with $A_1$ as the target matrix. Of course, $S(\widehat{N_k})$ is a dense matrix that is too expensive to compute in practice, but our purpose for this experiment is to investigate any underlying structure to the magnitude of the elements, particularly those nearer to the main diagonal. We apply a global sparsification strategy \cite{chow2000priori} to reveal an interesting behavior in the resulting sparsity patterns; see Figure \ref{THT}. In particular, we observe an increase in the number of nonzeros in the sparsity pattern of the sparsified exact map, demonstrating a decay-like behavior away from the main diagonal as the global threshold permits larger magnitude elements in the sparsified exact map. 

In fact, this is derived from the fact that although the inverse of a matrix is very dense, the matrix entries decay as they move further from the diagonal, especially for matrices representing the discrete Green's function \cite{chow2000priori,huckle1999approximate,tang1999toward} associated with the Laplace operator. However, the work in \cite{carr2021preconditioning} mainly focused on choosing sparsity patterns to maximize the efficiency of the map in sequential environments, with an emphasis on balancing the cost of computing and applying the map with the increased cost of GMRES iterations when solving the preconditioned systems using (potentially) dense maps. Indeed, in \cite{carr2021preconditioning, GriDesGug2020}, it is shown that $S(A_0)$ is sufficient to efficiently compute an accurate enough map for sequences of matrices from several applications, including THT. In this paper, we focus on identifying sparsity patterns that lead to more accurate approximate maps to drive future work that will focus on efficient approaches for computing those sparsity patterns identified in this paper, particularly in high-performance and parallel computing environments. Rather than sparsifying the exact map, which is far too computationally expensive, our THT example motivates the use of higher powers (level 1 or level 2 neighbors in $G(A_k)$) of the sparsified version of $A_k$\cite{chow2000priori}, denoted $\widetilde{A_k}$, as the sparsity pattern for SAM.

We apply three strategies for sparsifying $A_k$ to obtain $\widetilde{A_k}$ using different global and local threshold parameters inspired by previous works\cite{chow2000priori,chow2001parallel}. First, we use a global threshold, $thresh$, to construct $\Tilde{A_k}$ as\cite{chow2001parallel}:
\begin{equation}
    \widetilde{A_{ij}} = \begin{cases}
        1 & \text{if $i=j$ or $\mid (D^{\frac{1}{2}}AD^{\frac{1}{2}})_{ij} \mid > thresh$} \\ 
        0 & \text{otherwise}
    \end{cases}
\end{equation}
where the diagonal matrix $D$ is
\begin{equation}
    D_{ii} = \begin{cases}
        \mid A_{ii} \mid & \text{if $\mid A_{ii} \mid > 0$} \\
        1 & \text{otherwise.}
    \end{cases}
\end{equation}
However, while our main focus is on the quality of the map, we are still cognizant that these maps should not become too dense. To achieve good performance, we need to carefully consider the magnitude of the elements of the underlying system when choosing a global threshold parameter. More specifically, a global threshold that is too small will not impose enough sparsification to justify the use of the SAM. We refer to this strategy in our results as ``global {\it thresh}."

The second strategy imposes a local/column threshold, $\tau$, based on the magnitude of the largest matrix entry in that column. This can be thought of as the percentage of entries that should be allowed in each column of the sparsity pattern based on the largest absolute value entry in the corresponding column. The formula for sparsifying $A_k$ can be defined as\cite{chow2000priori},
\begin{equation}
    \widetilde{A_{ij}} = \begin{cases}
        1 & \text{if $i=j$ or $|A_{ij}| > (1 - \tau) \max_i |A_{ij}| $, where $0 \le \tau \le 1$} \\ 
        0 & \text{otherwise.}
    \end{cases}
\end{equation}
Because of the local nature of this dropping strategy, we conjecture that the column threshold approach should better capture the underlying system's behavior and produce effective sparsity patterns for the SAMs. Additionally, the sparsification for each column is independent of other columns, which makes it very suitable for implementations in parallel computing environments. This will be explored as part of future work. This will be explored as part of future work. We refer to this strategy in our results as ``col threshold, $\tau$."

The final strategy specifies how many nonzeros we want in each column of the sparsified matrix $\widetilde{A_k}$. Here, we define a parameter $lfil$ that limits the (largest in magnitude) entries of the column \cite{chow2000priori}. However, if the $lfil$ parameter is greater than the actual number of nonzeros in the column, it does not introduce new nonzeros. The formulation with the $lfil$ parameter is given by,
\begin{equation}
    \widetilde{A_{ij}} = \begin{cases}
        1 & \text{if $i=j$ or $|A_{ij}| \ge lfil(\max_i |A_{ij}|) $} \\ 
        0 & \text{otherwise.}
    \end{cases}
\end{equation}
We refer to this strategy in our results as ``fixed nz {\it lfil}."

\section{Numerical Experiments}

We apply the sparsity patterns with three different dropping strategies introduced above to two different applications that generate a sequence of linear systems: the numerical solution of the nonlinear convection-diffusion equation and discretized Euler equations. We describe these applications in section \ref{sec:cd2d} and \ref{sec:NWP}, respectively. We also compare against the sparsity pattern used in \cite{carr2021preconditioning}, $S(A_0)$. As our primary metric for discussing the quality of the map, we examine the Frobenius norm of the relative residual, $\|R_k\|_F/\|A_0\|_F = \|A_kN_k - A_0\|_F/\|A_0\|_F$, for each choice of sparsity pattern. We conduct our experiments in MATLAB (version R2022b) and do not consider the execution time for the computation. Rather, we focus on the quality of the sparsity patterns (i.e., large or small residual norm), and note that future work will explore efficient implementations for determining sparsity patterns in high performance computing environments. We also do not explicitly compare the performance of SAM as a preconditioner update for solving the sequence of linear systems. Empirical studies in \cite{carr2021preconditioning} and theoretical considerations of the residual norm in \cite{GriDesGug2020} indicate a direct correlation between the size of the residual norm (and the rate of increase/decrease across the sequence) and the effectiveness of SAMs as a preconditioner update with respect to total iterations to solve the entire sequence of linear systems. So, where our analysis will determine effective sparsity patterns that reduce $\|R_k\|_F/\|A_0\|_F$, we expect that those with the smallest residual norm will also result in better preconditioner updates. Future work will seek to verify and characterize this relationship further. 

In our experiments discussed next, we use a range of values for global threshold $thresh = 0.1$ down to $0.0001$ with decreasing orders of magnitude and let the local threshold s $\tau = 0.6, 0.7$, and $0.8$. For the fixed non-zero sparsity pattern, we let $lfil = 3$ and $5$ for the CD2D application and $lfil = 5$ and $7$ for the NWP application. Then, when applying the sparsity patterns, we first sparsify the source matrix $A_k$ to obtain $\widetilde{A_K}$, take the level 1 neighbor of $\widetilde{A_k}$, and compute the SAM using each of the three sparsity patterns. Note again that we map each back to the target matrix $A_0$. 

\subsection{Two-Dimensional Nonlinear Convection-Diffusion}\label{sec:cd2d}

The first application we chose is the nonlinear two-dimensional convection-diffusion equation as described in \cite{CarDiss2021}
 \begin{equation}\label{eq:convDif}
   -\nabla \cdot ((\eta + \gamma u^2)\nabla u) + r u_x + s u_y +t u = f,
 \end{equation} 
with $\eta = 0.1$, $\gamma = 1$, $r = s= 1$, and $t = f = 0$.  On the domain $[0,1]\times [0,1]$, we use the boundary conditions 
\begin{equation*}
u(x,1) = 0, ~~~~ u(x,0) = 0, ~~~~ u(1,y) = 0, ~~~~ u(0,y) = 0.2+y(1-y^2).
\end{equation*} 
We use second order central finite differences to discretize (\ref{eq:convDif}); the MATLAB code for the discretization scheme is provided in \cite[Appendix B.1]{CarDiss2021}. Since the PDE is nonlinear in the diffusion coefficient, Newton's method with line search is used to compute 
\begin{equation*}
    u^{(k+1)} = u^{(k)} - \alpha_{k+1}J(u^{(k)})^{-1}F(u^{(k)}),
 \end{equation*}
resulting in a sequence of linear systems of the form $F(u^{(k)}) = \nabla f(u^{(k)})$ where $u_{(k+1)}$ is the approximate solution in the $(k+1)^{st}$ Newton step, $J(u^{(k)})$ is the Jacobian evaluated at $u^{(k)}$, and $F(u^{(k)}) = \nabla f(u^{(k)})$. We refer to this as the CD2D application in our results and discussions.

\begin{figure}[hh]
\begin{center}
    \captionsetup{justification=centering}
    \includegraphics[width=13cm]{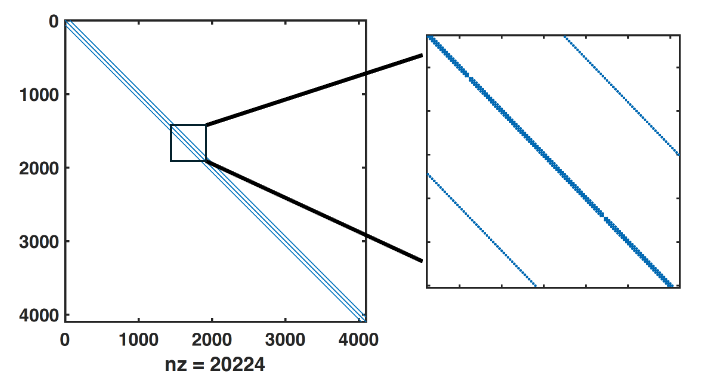}    
    \caption{Sparsity patterns for the CD2D matrices.}
    \label{fig:CD2Dpatt}
\end{center}
\end{figure}

For our experiments, we consider 66 matrices with $n = 4096$ satisfying (\ref{property}). The matrices $F(u^{(k)})$ exhibit the block tridiagonal sparsity pattern typical of the five-point stencil used for this discretization; see figure \ref{fig:CD2Dpatt}.

\begin{figure}[h]
\begin{minipage}{160mm}
    \captionsetup{justification=centering}
    \includegraphics[width=\linewidth,height=90mm]{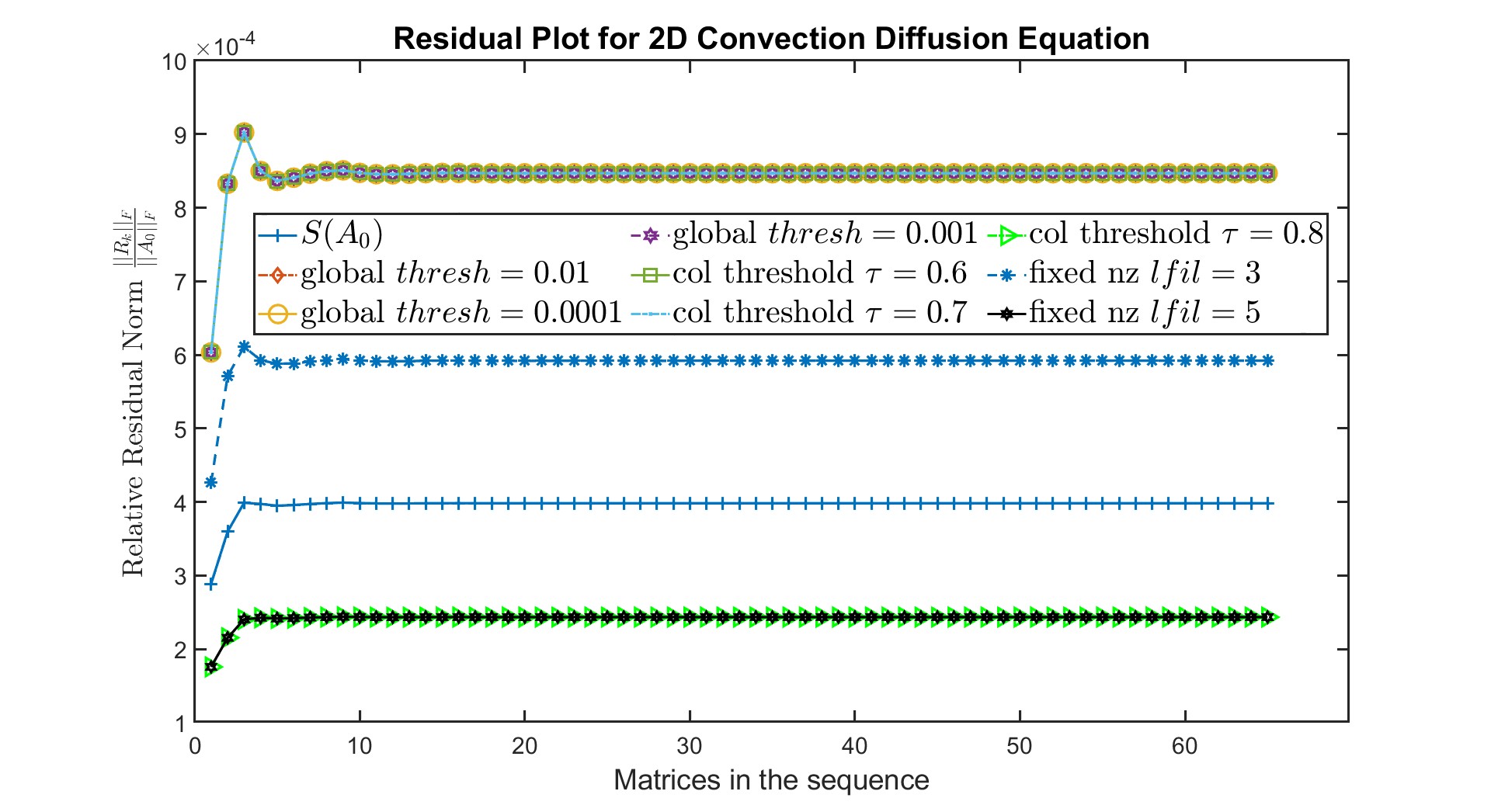}
    \caption{Relative residual $\frac{||R_k||_F}{||A_0||_F}$ in CD2D problem. Approximate maps $N_k$ are computed using the square of $S(\widetilde{A_k})$.}
    \label{cd2d_relative_resdual}
    \end{minipage}
\hfil
\begin{minipage}{160mm}
    \captionsetup{justification=centering}
    \includegraphics[width=\linewidth,height=90mm]{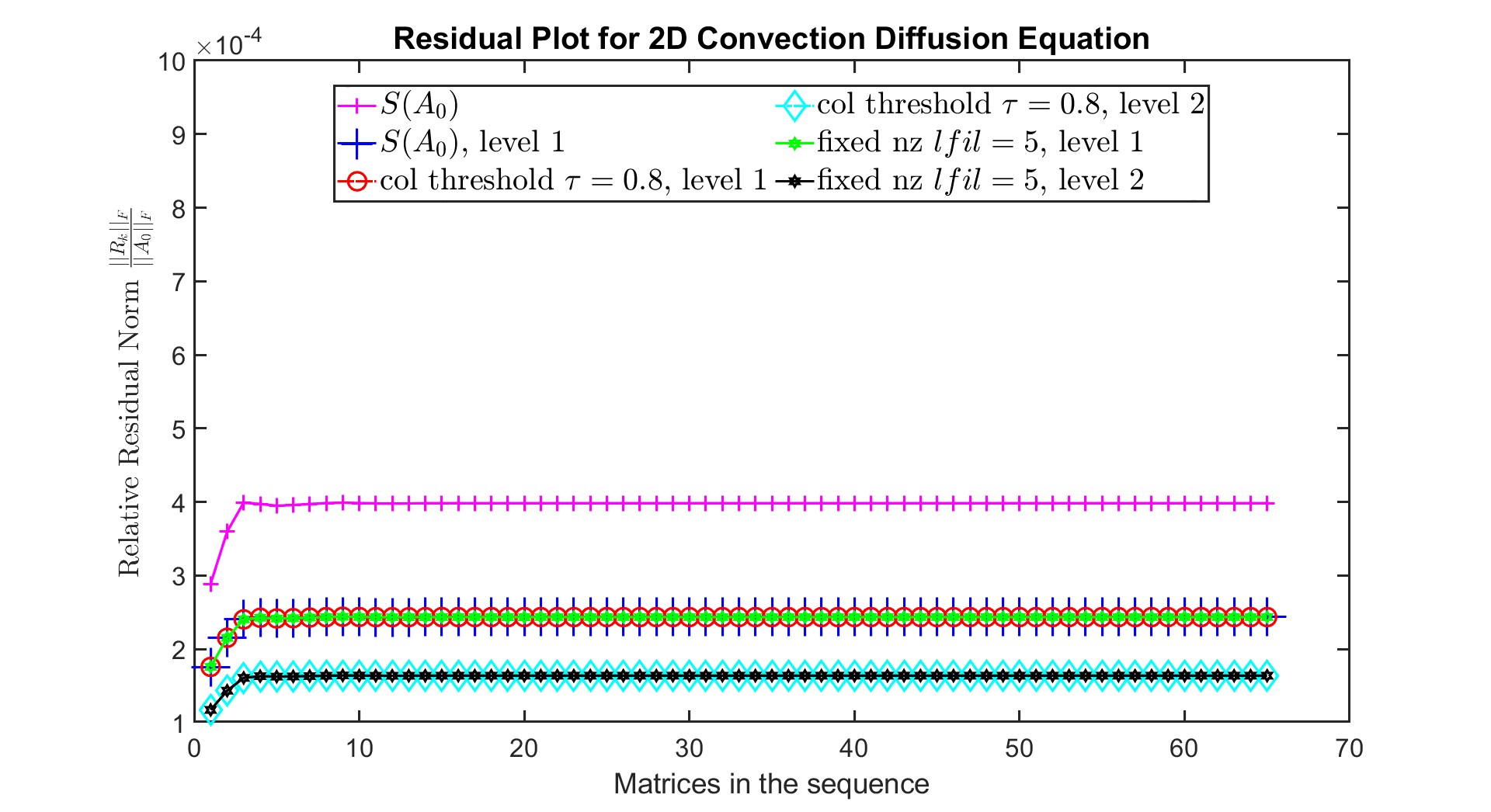}
    \caption{Comparison of relative residual norms between level 1 level 2 neighbors of $G(\widetilde{A_k})$.}
    \label{cd2d_relative_resdua_level_2}
\end{minipage}
\end{figure}

In figure \ref{cd2d_relative_resdual}, we plot the Frobenius norm of the relative residual $\frac{||R_k||_F}{||A_0||_F}$, where $||R_k|| = ||A_k N_k - A_0||_F$, after computing $N_k$ using these sparsity patterns in figure \ref{cd2d_relative_resdual}. Additionally, the number of non-zeros in the sparsity pattern of the computed SAMs is reported in table \ref{cd2d_nnz} for a few selected matrices. First, we look at the relative residuals for $N_k$ computed with level 1 neighbors of $\widetilde{A_k}$ sparsified by the global threshold parameter $thresh = 0.01, 0.001, 0.0001$ in figure \ref{cd2d_relative_resdual}. Even though the values of $thresh$ were different, they produced identical relative residuals, all of which are higher than the relative residual obtained when using the sparsity pattern of the target matrix $S(A_0)$ for computing $N_k$. This is expected because, as seen in table \ref{cd2d_nnz}, all the chosen parameter values for global $thresh$, have the same number of nonzeros in the sparsity pattern of the approximate map, which are comparatively fewer than the number of non-zeros when using $S(A_0)$ for computing the maps. Interestingly, the number of nonzeros does not change significantly for the different approximate maps in this sequence. This suggests that the discrepancy between the large magnitude entries and the small magnitude entries is fairly significant. As a result, we could not capture the smaller entries with $thresh$ value up to $0.0001$. For our experiments, we did not go beyond $thresh = 0.0001$; however, smaller values of $thresh$ might be able to capture more non-zero entries from the matrices.

Next, we compute the SAMs by imposing a sparsity pattern of level 1 neighbors of $\widetilde{A_k}$ sparsified by column threshold $\tau = 0.6, 0.7$, and $0.8$. For $tau = 0.6$ and $0.7$, the relative residuals are the same but also higher than those for $N_k$ computed with $S(A_0)$, as seen in figure \ref{cd2d_relative_resdual}. A similar pattern is observed in Table \ref{cd2d_nnz}, where the number of non-zero entries remains the same and does not change between matrices in the sequence for these two values of $\tau$. On the other hand, for $\tau = 0.8$, the relative residual for $N_k$ with $S(\widetilde{A_k})$ is smaller than the relative residual for $N_k$ computed with $S(A_0)$. Thus, a column threshold $\tau = 0.8$ is a good choice for sparsifying the source matrices and thus acts as a good sparsity pattern for computing the approximate maps. Interestingly, even for $\tau = 0.8$, the number of non-zeros for matrices in the sequence does not change. This gives us insight into the patterns of the matrices generated by the CD2D problem. In addition to the sparsity patterns of the matrices being identical, the magnitude of the entries in the matrices in the sequence does not change significantly. This allows us to conduct a lower bound analysis for the magnitude of entries to satisfy \ref{property} in future works.

Finally, we apply fixed $lfil = 3$ and $lfil = 5$ and impose the sparsity pattern of the resulting $\widetilde{A_k}$ to compute the approximate maps. Observing figure \ref{cd2d_nnz}, for $lfil = 3$, the relative residuals are larger than those of $N_k$ computed with $S(A_0)$, although they are not as high as those for the global column threshold choices discussed earlier. For $lfil = 5$, relative residuals are indeed smaller than those for approximate maps computed with the sparsity pattern of the target matrix. However, we must consider that our CD2D problem uses a five-point stencil for discretization. As a result, each column will contain at most five non-zero entries. So, when we take the level 1 neighbor of $S(\widetilde{A_k})$ with $lfil = 5$, we are effectively taking the level 1 neighbor of the sparsity pattern of the source matrices, which is the same as the target matrix. When we observe Table \ref{cd2d_nnz}, the number of non-zeros with $lfil = 5$ does not change for matrices in the sequence, whereas for $lfil = 3$, that number changes. Additionally, comparing the number of non-zeros of the maps for column threshold $\tau = 0.8$ with other choices in Table \ref{cd2d_nnz}, the numbers are identical to those with $lfil = 5$. This suggests they generate the same sparsity patterns for $\widetilde{A_k}$, which is the level 1 neighbor of all the non-zero entries in the source matrices. This means the magnitude of the off-diagonal entries of the matrices is within $20\%$ of the largest matrix entries, i.e., the main diagonal, but not within $30\%$ since they were not selected for the sparsity pattern when $\tau = 0.7$.

\begin{table}
\centering
\captionsetup{justification=centering}
\caption{Number of non-zeros in the sparsity patterns of the map, $S(N_k)$ for level 1 neighbors of different sparsification choices applied to the CD2D problem. Only a selected few matrices in the sequence are shown.}
\label{cd2d_nnz}
\begin{tabular}{@{}lllllll@{}}
\hline
Sparsity Pattern & $nnz(N_2)$ & $nnz(N_{5})$ & $nnz(N_{10})$ & $nnz(N_{30})$ & $nnz(N_{45})$ & $nnz(N_{66})$ \\ \hline
$S(A_0)$ & $20224$ & $20224$ & $20224$ & $20224$ & $20224$ & $20224$ \\
global $thresh = 0.01$ & $4103$ & $4157$ & $4142$ & $4098$ & $4105$ & $4106$ \\
global $thresh = 0.001$ & $4103$ & $4157$ & $4142$ & $4098$ & $4105$ & $4106$ \\ 
global $thresh = 0.0001$ & $4103$ & $4157$ & $4142$ & $4098$ & $4105$ & $4106$ \\
col threshold $tau = 0.6$ & $4096$ & $4096$ & $4096$ & $4096$ & $4096$ & $4096$ \\
col threshold $tau = 0.7$ & $4096$ & $4096$ & $4096$ & $4096$ & $4096$ & $4096$ \\
col threshold $tau = 0.8$ & $51972$ & $51972$ & $51972$ & $51972$ & $51972$ & $51972$ \\
fixed nz $lfil = 3$ & $20662$ & $20781$ & $20758$ & $20672$ & $20673$ & $20673$ \\
fixed nz $lfil = 5$ & $51972$ & $51972$ & $51972$ & $51972$ & $51972$ & $51972$ \\ \hline
\end{tabular}
\end{table}

We further investigated the relative residual norms with level 2 neighbors of $G(\widetilde{A_k})$ in figure \ref{cd2d_relative_resdua_level_2} and reported the number of non-zeros in $S(N_k)$ in table \ref{cd2d_nnz_level2}. For this analysis, we chose only the column threshold $\tau = 0.8$ and fixed nonzero $lfil = 5$ for computing $\widetilde{A_k}$ since the relative residuals for approximate maps with these choices were smaller than the approximate maps with the sparsity pattern of the target matrix. As expected, taking higher powers of the sparsity pattern, i.e., taking higher level neighbors of $S(\widetilde{A_k})$, generates SAMs with better approximation. Sparsity patterns with level 2 neighbors generate better approximations compared to the level 1 neighbor versions of all the corresponding choices, as seen in Figure \ref{cd2d_relative_resdua_level_2}. The relative residual of the level 2 neighbor for column threshold $\tau = 0.8$ and $lfil = 5$ is smaller than both the relative residuals from the sparsity patterns of the target matrix and its level 1 neighbor. Additionally, the values of the relative residuals are the same for all the computed SAMs, exhibiting the same behavior as their level 1 versions. They also have the same number of non-zeros in the sparsity pattern of the SAMs, as seen in Table \ref{cd2d_nnz_level2}. Furthermore, the level 1 neighbor of the target matrix itself has an identical relative residual compared to the level 1 neighbor with column threshold $tau = 0.8$ and $lfil = 3$ in figure \ref{cd2d_nnz_level2}. This pattern is further reflected in table \ref{cd2d_nnz_level2}, where all strategies result in the same number of non-zeros for all the maps. This confirms our assumption that the sparsity pattern generated by the sparsification techniques with column threshold $\tau = 0.8$ and $lfil = 5$ is essentially the sparsity pattern of the original matrices. They could not sparsify the matrices adequately.

\begin{table}
\centering
\captionsetup{justification=centering}
\caption{Number of non-zeros in the sparsity patterns of the map, $S(N_k)$ for level 1 and level 2 neighbors of selected sparsification choices applied to the CD2D problem. Only a selected few matrices in the sequence are shown.}
\label{cd2d_nnz_level2}
\begin{tabular}{@{}lllllll@{}}
\hline
Sparsity Pattern & $nnz(N_2)$ & $nnz(N_{5})$ & $nnz(N_{10})$ & $nnz(N_{30})$ & $nnz(N_{45})$ & $nnz(N_{66})$ \\ \hline
$S(A_0)$ & $20224$ & $20224$ & $20224$ & $20224$ & $20224$ & $20224$ \\
$S(A_0)$, level 1 & $51972$ & $51972$ & $51972$ & $51972$ & $51972$ & $51972$ \\
col threshold $tau = 0.8$, level 1 & $51972$ & $51972$ & $51972$ & $51972$ & $51972$ & $51972$ \\
col threshold $tau = 0.8$, level 2 & $98836$ & $98836$ & $98836$ & $98836$ & $98836$ & $98836$ \\
fixed nz $lfil = 5$, level 1 & $51972$ & $51972$ & $51972$ & $51972$ & $51972$ & $51972$ \\
fixed nz $lfil = 5$, level 2 & $98836$ & $98836$ & $98836$ & $98836$ & $98836$ & $98836$ \\ \hline
\end{tabular}
\end{table}

CD2D linear systems exhibit very simple structures in terms of both sparsity patterns and the magnitude of matrix entries. The changes in the magnitudes of matrix entries between different matrices in the sequence are also minimal. The sparsity patterns remain consistent across matrices in the sequence. Due to the second-order central finite difference discretization with a five-point stencil, the matrices are very sparse, containing only five non-zero entries per column. Owing to the diagonal dominance property of CD2D problems \cite{Golub1996}, the largest matrix entries are along or near the main diagonal, with the distant off-diagonal entries having significantly smaller magnitudes. Given these constraints, it is difficult to sparsify the matrices further, resulting in mostly identical sparsity pattern choices, as evident from Figures \ref{cd2d_relative_resdual} and \ref{cd2d_relative_resdua_level_2}, and Tables \ref{cd2d_nnz} and \ref{cd2d_nnz_level2}, as well as our previous discussion. In these types of model problems, there is little difference between using the sparsity pattern of the target matrix or the source matrix. However, to achieve higher accuracy in the approximate maps, utilizing the level 1 and level 2 neighbors of these source matrices, $S(A_k)$, is effective.

\subsection{Euler Equations for Warm Bubble Rising}\label{sec:NWP}

\begin{figure}[h]
\captionsetup{justification=centering}
  \centering
    \includegraphics[width=0.8\textwidth]{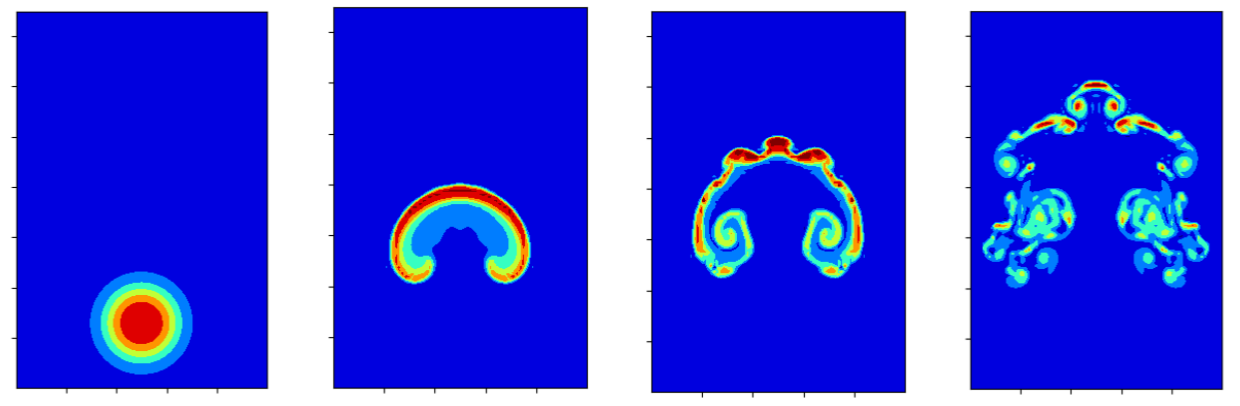}
\caption{Warm bubble rising in the atmosphere.}\label{fig:BubbleVisual}
\end{figure}

\begin{figure}[hh]
\begin{center}
    \captionsetup{justification=centering}
    \begin{subfigure}{.33\textwidth}
        \includegraphics[width=\linewidth]{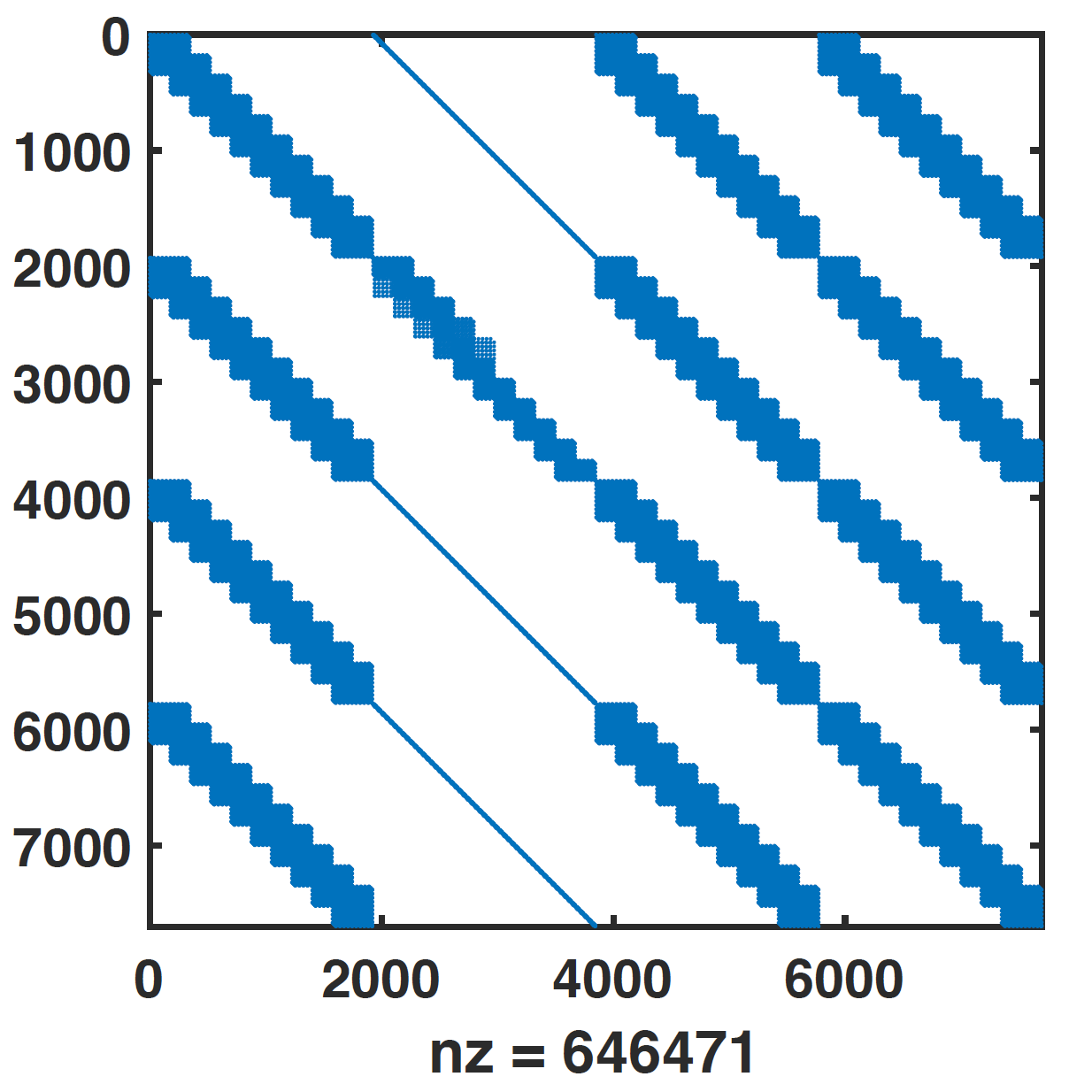}
		\subcaption{$S(A_0)$}\label{bubble_A0}
    \end{subfigure}        
    \begin{subfigure}{.33\textwidth}
        \includegraphics[width=\linewidth]{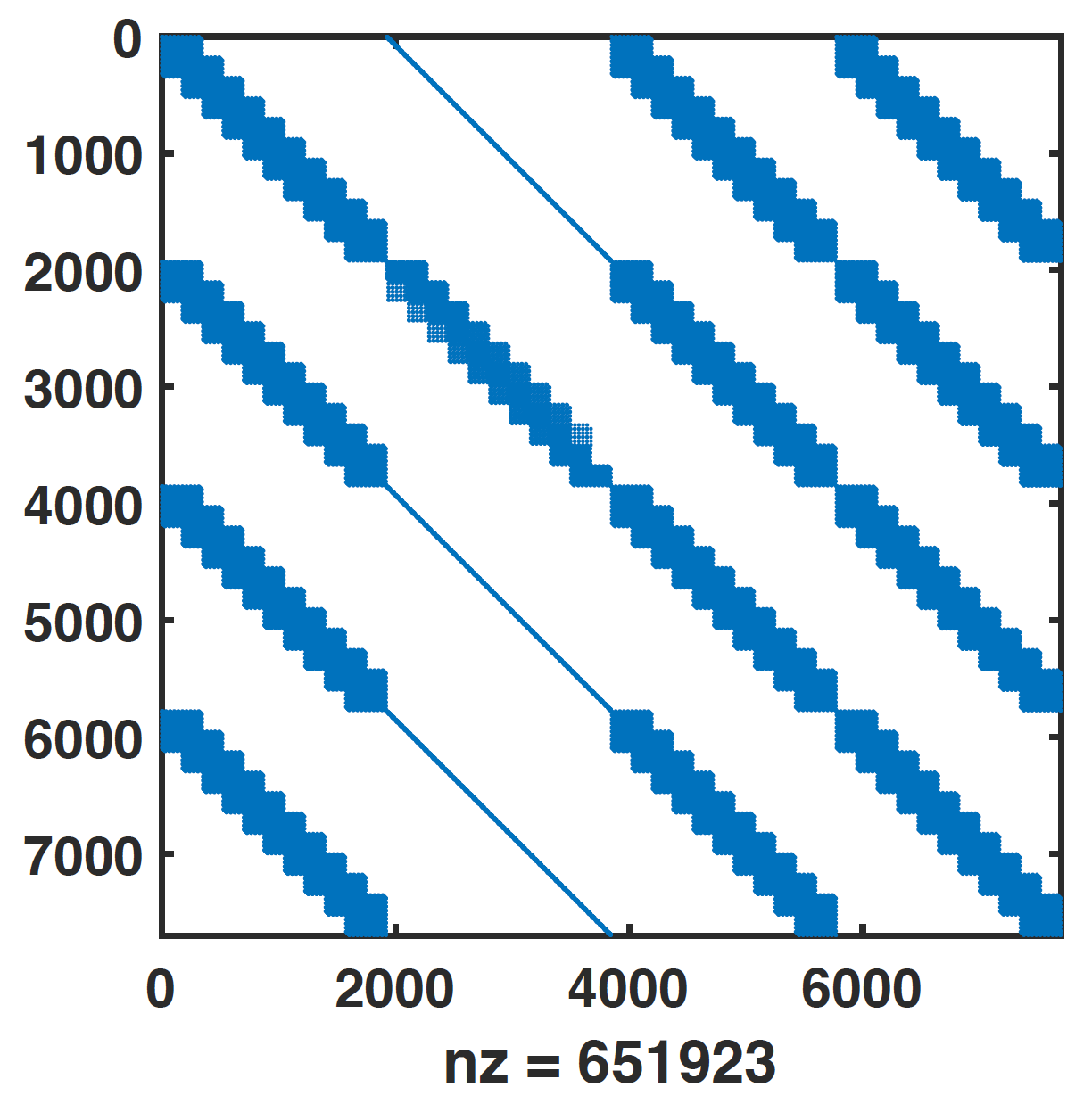}
		\subcaption{$S(A_{11})$}\label{bubble_A10}
    \end{subfigure}
    \begin{subfigure}{.33\textwidth}
        \includegraphics[width=\linewidth]{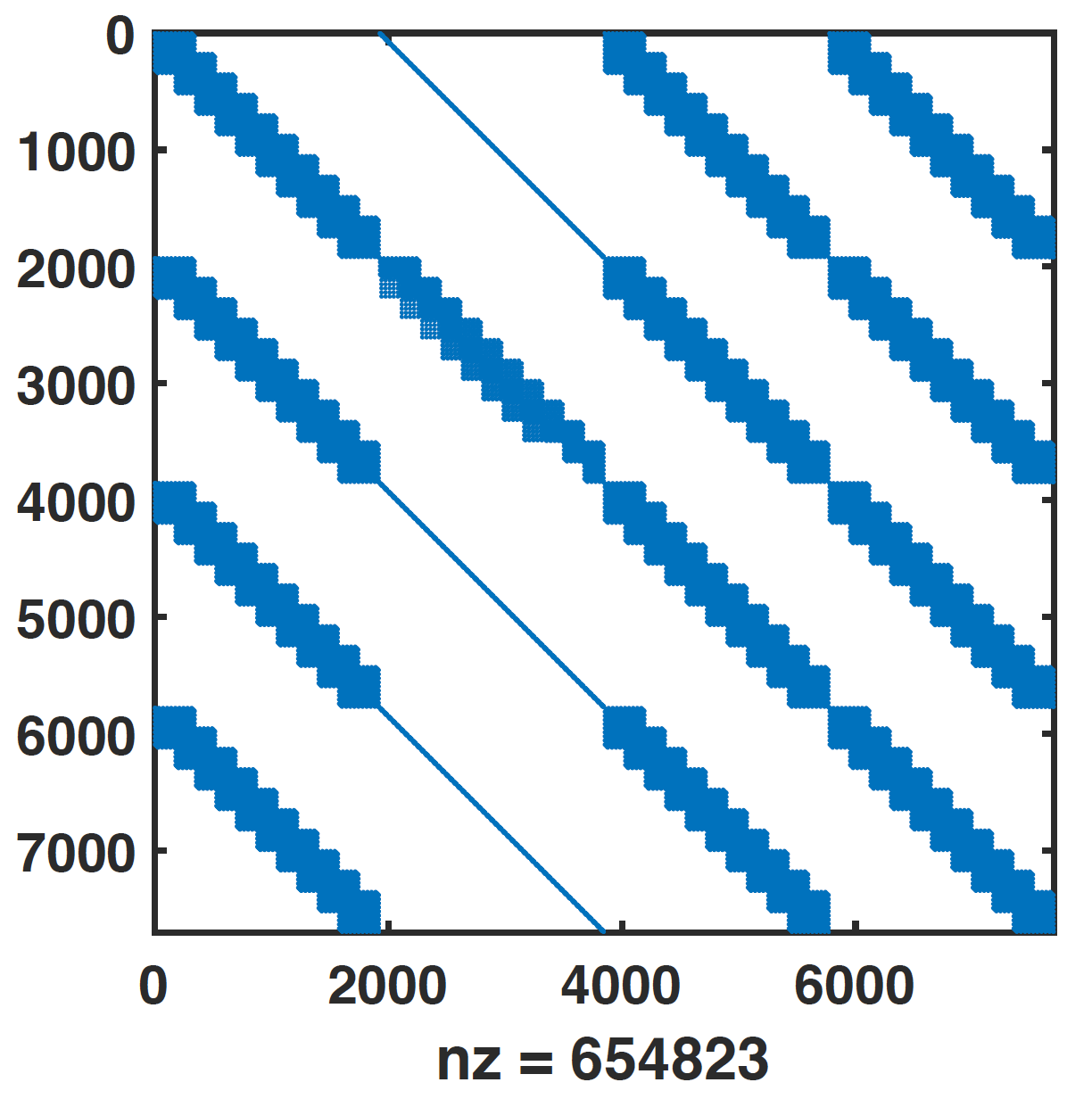}
		\subcaption{$S(A_{19})$}\label{bubble_A19}
    \end{subfigure}
  \caption{Sparsity patterns for the NWP matrices.}
  \label{fig:bubblepatt}
\end{center}
\end{figure}

The second application stems from an implicit time discretization of the Euler equations on a rotating sphere\cite{gaudreault2022high} and is frequently used as an initial step in the design of a numerical weather prediction (NWP)\cite{NWP}. We use a simplified model based on the Euler equations that simulate the dynamics of a warm bubble rising in a dry isentropic (i.e., constant entropy) atmosphere; see figure \ref{fig:BubbleVisual} for a visualization of the simulation and \cite{robert1993bubble} for a detailed description of the problem. 
This underlying PDE is discretized using the finite element method on a two-dimensional Cartesian grid, and a second-order Rosenbrock scheme \cite[Eq. 8]{dallerit2023second} is run using code available from \cite{NWPgithub}. This is done for 216 steps with a timestep of 5 seconds and for each step, a linear system must be solved. For our experiments, we consider the first 20 coefficient matrices in this sequence, with $n = 7700$, and refer to this application as the NWP application. The sparsity pattern of these matrices is much denser than the other applications considered in this paper and they do not necessarily respect (\ref{property}). Figure \ref{fig:bubblepatt} shows the general sparsity pattern for the matrices in this sequence, noting that while each matrix in this sequence more or less maintains this same pattern, some matrices later in the sequence may not have nonzeros in certain locations, whereas earlier matrices in the sequence do, see figures \ref{bubble_A0}, \ref{bubble_A10}, and \ref{bubble_A19}. In particular, we highlight the main block diagonal of each matrix shown. Note that $S(A_{11})$ has a denser main diagonal block than $S(A_{19})$ and so $S(A_{11}) \not \subseteq S(A_{19})$; this is the case for many matrices along the sequence.  Nevertheless, we include this application to observe the effect of the sparsification strategies investigated in this paper as they are applied to matrices whose sparsity patterns do not satisfy the subset property (\ref{property}).

\begin{figure}
\centering
\captionsetup{justification=centering}
\includegraphics[width=16cm, height=9cm]{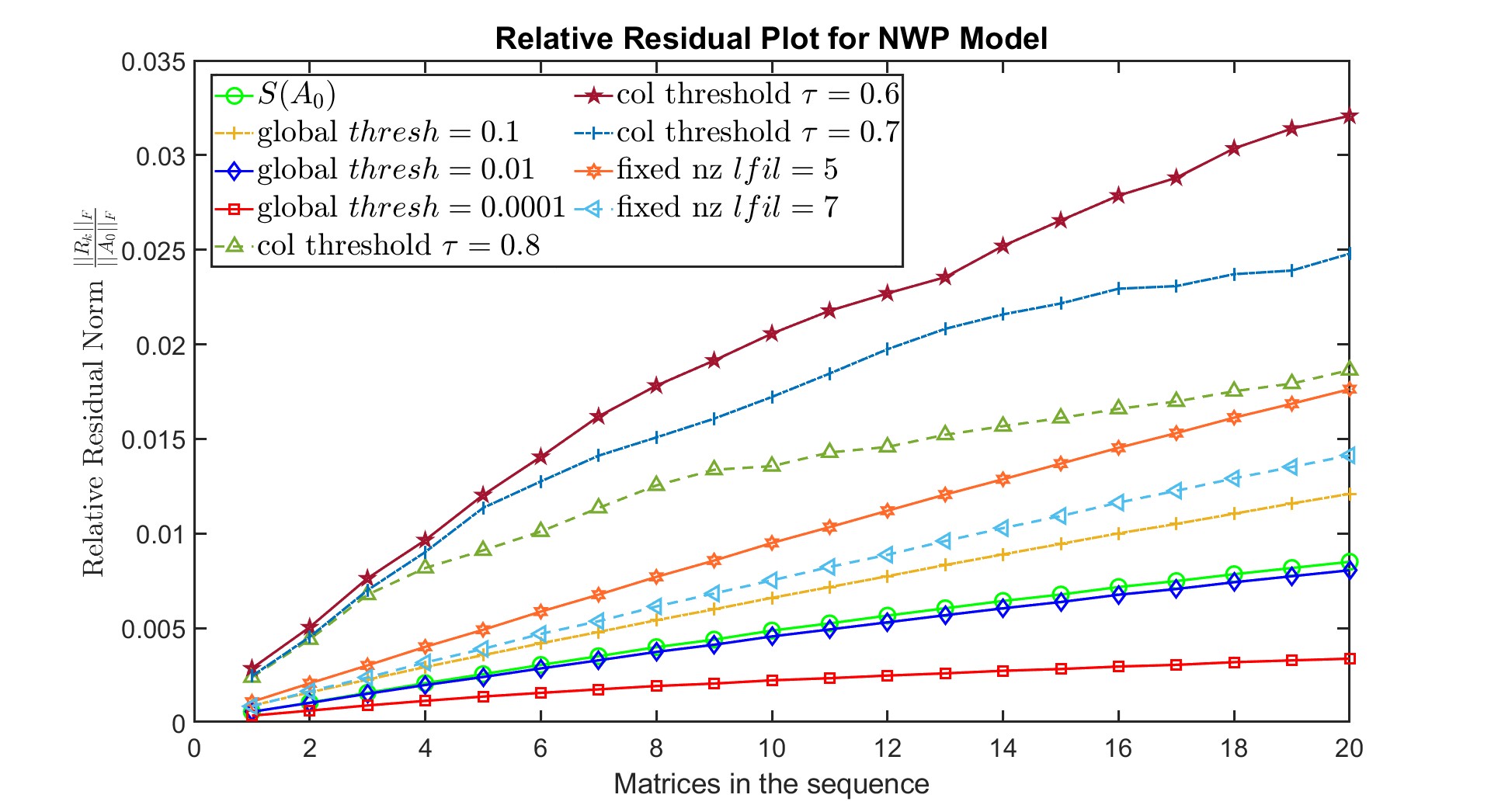}
\caption{Relative residual $\frac{||R_k||_F}{||A_0||_F}$ for NWP model. Approximate maps $N_k$ are computed using the square of $S(\widetilde{A_k})$.}
\label{bubble_map_relative_residual}
\end{figure}

Similar to the CD2D problem, we compute $N_k$ with different sparsified matrix $\widetilde{A_k}$ and plot the relative residual norm in figure \ref{bubble_map_relative_residual} and report the number of non-zeros in the map in table \ref{nwp_nnz}. Sparsifying $A_k$ with a global threshold value $thresh = 0.1, 0.01$ and $0.0001$ and computing the SAMs with the resulting $S(\widetilde{A_k})$ produces fairly small residuals compared with the SAMs computed using the initial sparsity pattern $S(A_0)$. For $thresh = 0.01$ and $0.0001$, the relative residual norm was smaller than using the sparsity pattern of the original matrices, $S(A_0)$. The number of non-zeros in table \ref{nwp_nnz} for these global sparsification choices are different for different values. This is expected since the sparsity pattern and magnitude of the NWP matrices are changing for each matrix in the sequence even though they do not satisfy \ref{property}. In this application, a decreasing value of $thresh$ resulted in an increased number of non-zeros in the sparsity pattern. Interestingly, we observe that the number of non-zeros also demonstrates an increasing trend across the matrices in the sequence for all values of $thresh$. This suggests that with each matrix in the sequence, the magnitude of the matrix entries is increasing in general.

\begin{table}
\centering
\captionsetup{justification=centering}
\caption{Number of non-zeros in the sparsity patterns of the map, $S(N_k)$ for level 1 neighbors of different sparsification choices applied to the NWP problem. Only a selected few matrices in the sequence are shown.}
\label{nwp_nnz}
\begin{tabular}{@{}lllllll@{}}
\hline
Sparsity Pattern & $nnz(N_2)$ & $nnz(N_{5})$ & $nnz(N_{10})$ & $nnz(N_{15})$ & $nnz(N_{20})$ & $nnz(N_{21})$ \\ \hline
$S(A_0)$ & $646471$ & $646471$ & $646471$ & $646471$ & $646471$ & $646471$ \\
global $thresh = 0.1$ & $438487$ & $463415$ & $483151$ & $496579$ & $506728$ & $508110$ \\
global $thresh = 0.01$ & $684883$ & $735602$ & $775809$ & $813448$ & $839414$ & $834187$ \\ 
global $thresh = 0.0001$ & $1389705$ & $1573152$ & $1747478$ & $1870605$ & $1963174$ & $1966994$ \\
col threshold $tau = 0.6$ & $95757$ & $95812$ & $96283$ & $97071$ & $97877$ & $97990$ \\
col threshold $tau = 0.7$ & $130153$ & $130314$ & $131525$ & $132959$ & $134511$ & $134814$ \\
col threshold $tau = 0.8$ & $193726$ & $194469$ & $197155$ & $200554$ & $203806$ & $204420$ \\
fixed nz $lfil = 5$ & $146310$ & $146276$ & $147831$ & $148935$ & $149066$ & $148312$ \\
fixed nz $lfil = 7$ & $240876$ & $239481$ & $241456$ & $243120$ & $243589$ & $242075$ \\ \hline
\end{tabular}
\end{table}

Next, we computed SAMs with $\widetilde{A_k}$ sparsified by column threshold $\tau = 0.6, 0.7$, and $0.8$ that generated larger residuals than using $S(A_0)$ as a sparsity pattern for the map. For $\tau = 0.8$, the relative residual decreases quite significantly for matrices that are later in the sequence compared to $\tau = 0.6$ and $0.7$. This can happen because of capturing more non-zero entries that allow for a more accurate approximation of the exact map. Referring to table \ref{nwp_nnz}, they maintain the same trend of increasing the number of non-zeros in the sparsity pattern of the map with each matrix in the sequence.

Finally, we observe the relative residual plots for fixed nonzero sparsification of the source matrix with parameters $lfil = 5$ and $7$. Both produce larger relative residuals compared to using $S(A_0)$ to compute SAMs but they are not significantly higher.  This is especially interesting considering the average number of nonzeros per column in an NWP matrix can be around 85. This suggests an aggressive dropping strategy to determine optimal sparsity patterns for this application can likely be used with good results. In this case, the cost to compute and apply a SAM with 5 or 7 nonzeros per column is significantly less expensive than when imposing the sparsity pattern of the matrix itself, while resulting in only a slightly less accurate approximate map. The number of nonzeros for these choices in table \ref{nwp_nnz} are fairly consistent for different matrices in the sequence as we are now fixing the number of non-zeros each column can have. So, even though the magnitude of the matrix entries is increasing this sparsification technique might not be able to capture them.

NWP matrices contain a large number of non-zero entries. Even though they do not maintain the subset property, using a global threshold to sparsify the matrix captured the matrix entries with a larger magnitude. Level 1 neighbors of these matrix entries produce a sparsity pattern that produces a better approximation to exact maps than SAM's current implementation. This can be an area of future exploration to find a way to characterize the sparsity pattern for matrices that do maintain (\ref{property}).

\section{Conclusion}
Our objective is to analyze the sparsity pattern of the exact map $\widehat{N_k}$ through graphical representation and characterize a priori sparsity patterns for approximate maps $N_k$ that provide better approximations to the exact maps for sequences of matrices generated from discretized PDEs. Provided the matrices in the sequence follow (\ref{property}), we should use the sparsity pattern of the source matrix $A_k$ rather than the target matrix $A_0$. Furthermore, we applied three simple sparsification techniques for the source matrix $A_k$ based on \cite{chow2000priori} and used powers of the sparsified matrix $\widetilde{A_k}$ as the sparsity for SAMs. We applied these sparsity patterns to two applications, one that satisfies (\ref{property}) and one that does not. For both problems, several choices provide a better approximation than when using $S(A_0)$, although we could not converge on a single choice that works in general. This is very much problem-dependent, and as a future direction, we can look at specific discretization techniques applied to PDEs resulting in known patterns and potentially magnitudes of elements. Moreover, we will also look at more sophisticated ways of sparsifying the matrices that will capture the underlying behavior of the PDEs and provide a more efficient, but still accurate, approximation. Also, combining this with implementations in a parallel computing environment will help us compute more accurate and computationally efficient SAMs.

\def\bstname{pamm}
\bibliographystyle{pamm}
\bibliography{pamm-tpl}

\end{document}